\documentclass[epic,eepic,11pt]{amsart}

\usepackage{amsfonts,mathrsfs}
\usepackage[mathscr]{eucal}

\usepackage{amssymb}
\usepackage{amscd}
\usepackage{fancyhdr}

\usepackage{euler,eucal}

\DeclareMathAlphabet{\mathbf}{T1}{ppl}{bx}{n}
\DeclareMathAlphabet{\mathrm}{T1}{ppl}{m}{n}

\newcommand\note[1]%
{$^\dagger$\marginpar{\footnotesize{$^\dagger${#1}}}}

\def\({\left(}
\def\){\right)}
\def\<{\left<}
\def\>{\right>}

\def\newi {\sqrt{-1}\, }

\newtheorem{theorem}{Theorem}
\newtheorem{proposition}[theorem]{Proposition}
\newtheorem{lemma}[theorem]{Lemma}
\newtheorem{definition}[theorem]{Definition}

\theoremstyle{definition}
\newtheorem{example}[theorem]{Example}
\newtheorem{remark}[theorem]{Remark}

\newcommand\lie{\mathfrak}

\newcommand\g{\lie{g}}

\newcommand\C{\mathbb{C}}

\newcommand\CP{\mathbb{CP}}

\newcommand\J{\mathcal{J}}
\newcommand\cO{\mathcal{O}}

\newcommand     {\comment}[1]   {}
\newcommand{\mute}[2] {}
\newcommand     {\printname}[1] {}

\newcommand\func[1]{\operatorname{\mathrm{#1}}}
\newcommand\funclim[1]{\operatorname*{\mathrm{#1}}}

\newcommand\type{\func{type}}

\renewcommand\dim{\func{dim}}

\renewcommand\lim{\funclim{lim}}

\newcommand\sur{\mathrel{\to\kern-1.8ex\to}}
\newcommand\iso{\mathrel{\hookrightarrow\kern-1.8ex\to}}

\newcommand\longhookrightarrow{\lhook\joinrel\longrightarrow}

\newcommand\longsur{\mathrel{\longrightarrow\kern-1.8ex\to}}
\newcommand\longiso{\mathrel{\longhookrightarrow\kern-1.8ex\to}}

\begin{document}

\bibliographystyle{amsalpha}
\date{\today}

\title{Reduction of twisted generalized K\"ahler structure}

\author{Yi Lin, \,\,\,  Susan Tolman}
\address{Department of Mathematics, University
of Toronto, Canada, M5S2E4, yilin@math.toronto.edu}
\address{Department of Mathematics, University of Illinois at
Urbana-Champaign, Urbana, IL, USA, 61801, stolman@math.uiuc.edu}

\begin{abstract}
We show that our  earlier work in \cite{LT05}  extends to the
twisted case, that is, we defined a notion of moment map and
reduction in both twisted generalized complex geometry and twisted
generalized K\"ahler geometry.
\end{abstract}

\maketitle

This is a short note to show that our results in our previously
posted paper, {\em Symmetries in generalized K\"ahler geometry}
\cite{LT05}, can be easily extended to the case of twisted
generalized geometries with only minor modifications. In particular,
we define twisted generalized complex reduction and twisted
generalized K\"ahler reduction. We hope that this goes some way
towards providing the framework which Kapustin and Li suggested
would be useful \cite{KL}.

This note is not intended as a completely separate  paper, but
should rather be read in conjunction with \cite{LT05}. In
particular, we will not repeat the brief history of this subject
which can be found there, except to reiterate that it was first
introduced in \cite{H02} and developed much further in \cite{Gua04}.

However, it is important to mention here that  closely related works
which have  appeared since we first posted. We have not had an
opportunity to carefully read these papers but would like to at
least attempt to describe some basic similarities and differences.

Shengda Hu  wrote a related paper \cite{Hu05}, which was partially
inspired by an early version of our paper \cite{LT05} which we gave
him in early June. His paper includes a notion of twisted complex
reduction which is very similar to ours in the untwisted case. More
generally, he considers twisted complex structures in the framework
of Hamiltonian symmetry.  The work in the current note was completed
after his paper appeared, so also cannot be considered as fully
independent. Our main motivation is to demonstrate twisted
generalized K\"ahler reduction, which is not in \cite{Hu05}.
However, it is worth noting that even in the twisted generalized
complex case our results are slightly different.

In \cite{SX05}, Sti\'enon and Ping  independently develop  notions
of generalized complex and K\"ahler reduction which seems rather
different from ours in both perspective and techniques. In
particular, instead of working with generalized moment maps, they
consider quotients of arbitrary subsets; so our theorems  do not
appear in their paper (or vise-versa). They do not consider
bi-Hermitian or twisted structures.

In \cite{BCG05}, which is also independent of our work, Busztychn,
Cavalcanti, and Gualtieri define a quite general notion of the
reduction of Courant algebroids and Dirac structures which  includes
twisted generalized complex and twisted generalized K\"ahler
reduction as special cases.  They work with both arbitrary subsets
and moment maps, but again their techniques are different from ours.
They also construct two bi-Hermitian structures on $\CP^2$, but do
not construct  bi-Hermitian structures on the other spaces we
consider in \cite{LT05}.

We begin with a brief introduction, following \cite{Gua04}. (We
would like to mention that in addition to the many concrete
instances where we quote Gualtieri, our whole perspective on
generalized geometries was heavily influenced by his excellent
thesis on this subject.)

Let $M$ be a $n$ dimensional manifold. There is a natural metric of
type $(n,n)$ on $TM \oplus T^*M$ given by
$$\langle X + \alpha, Y + \beta \rangle = \frac{1}{2}
( \alpha(Y) + \beta(Y)),$$ which extends naturally to $T_\C M\oplus
T_\C ^*M=(TM\oplus T^*M)\otimes \C$.

A {\bf generalized almost complex structure} on a manifold $M$ is a
 bundle map $\mathcal{J} \colon TM \oplus T^*M \to TM \oplus T^*M$ which is orthogonal with respect to
 the natural metric defined above so that $\mathcal{J}^2 =
-1$ .

 For a closed form $H \in
\Omega^3(M)$, the {\bf $\bold{H}$-twisted Courant bracket} on $T_\C
M \oplus T_\C^* M$ is given by
$$ [X + \alpha, Y + \beta]_H  = [X,Y] + L_X \beta - L_Y \alpha
- \frac{1}{2} (d \iota_X \beta - d \iota_Y \alpha) + \iota_Y \iota_X
H. $$

A {\bf $\bold{H}$-twisted generalized complex  structure} on $M$ is
a generalized almost complex structure $\J$ so that the $\sqrt{-1}$
eigenbundle  of $\J$ is closed under the $H$-twisted Courant
bracket. And a {\bf $\bold{H}$-twisted generalized K\"{a}hler
structure} on a manifold $M$ is a pair of commuting
$\bold{H}$-twisted generalized complex structures $\J_1$ and $\J_2$
on $M$ so that $ - \J_1 \J_2$ is a positive definite bundle map on
$TM \oplus T^*M$ with respect to the natural metric.

As the following example shows, the concept of twisting is most
interesting when  $H$ represents a non-trivial cohomology class.

\begin{example}\cite{Gua04} Given a two-form $B$ on a manifold $M$, consider the  orthogonal
bundle map $TM \oplus T^*M \to TM \oplus T^*M$ defined by
$$e^B=\left( \begin{matrix} 1 & 0\\B & 1\end{matrix}\right),$$
where $B$ is regarded as a skew-symmetric map from $TM$ to $T^*M$.
If $\J$ is a generalized almost complex structure on $M$, then
$\mathcal{J}'= e^B \mathcal{J} e^{-B}$ is another generalized almost
complex structure on $M$, called the {\bf $\bold{B}$-transform} of
$\mathcal{J}$. For any  $B \in \Omega^2(M)$ and closed $H \in
\Omega^3(M)$, the $B$-transform of a $H$-twisted generalized complex
(or K\"ahler) structure is a $H  + dB$-twisted generalized complex
(or K\"ahler) structure.
\end{example}

It will also be convenient to have the following definition.

\begin{definition}
Let a compact Lie group $G$ act on a manifold $M$. The {\bf Cartan
model} for the equivariant cohomology of $M$ is defined as follows:
The degree $n$ co-chains are
$$ \Omega_G^{n}(M) = \bigoplus_i (\Omega^{n - 2i}(M)  \otimes S^i(\g^*))^G,$$
where $S^i$ denotes polynomials of degree $i$. The differential $d_G
:\Omega_G^n \to \Omega_G^{n+1}$ is defined by
$$d_G(\alpha \otimes p)(\xi) = (d \alpha  - \iota_{\xi_M} \alpha) p(\xi)
\qquad \mbox{for all} \ \xi \in \g,$$ where we think of
$\Omega_G^*(M)$ as the space of equivariant polynomial mappings from
$\g$ to $\Omega^*(M)$. (If $G$ acts on a vector space $A$, let $A^G$
denote the invariant subspace.) The {\bf equivariant cohomology} of
$M$ is $H_G^*(M) = H^*(\Omega_G^*,d_G)$.
\end{definition}

\begin{remark}
Let a compact Lie group $G$ act on a manifold $M$ so that it acts
freely on a submanifold $L \subset M$. There is a natural map,
called the  {\bf Kirwan map}
$$\kappa: H_G^*(M) \to H^*(L/G)$$  which
is the composition of the restriction map from $H_G^*(M)$ to
$H_G^*(L)$ with the natural isomorphism from $H_G^*(L)$ to
$H^*(L/G)$.

A form $B \in \Omega^n(M)$ is {\bf basic} if it is invariant and
$\iota_{\xi_M} B = 0$ for all $\xi \in \g$. Then $B$  {\bf descends}
to $\widetilde{B} \in \Omega^n(L/G)$, that is, the pull-back of
$\widetilde{B}$ to $L$ is the restriction of $B$.

If $B \in \Omega^n(M)^G \subset \Omega_G^n(M)$ is {\bf equivariantly
closed}, that is, $d_G B = 0$, then $B$ is closed and basic and
$\kappa[B] = [\widetilde{B}]$. More generally, if $\eta \in
\Omega_G^n(M)$ is equivariantly closed, there exists $\Gamma \in
\Omega_G^{n-1}(L)$ so that $\eta|_L + d_G \Gamma \in \Omega^n(L)^G
\subset \Omega_G^n(L)$. Since $\eta + d_G \Gamma$ is equivariantly
closed, it descends to $\widetilde{\eta} \in \Omega^n(L/G)$ and
$\kappa[\eta] = [\widetilde{\eta}]$.
\end{remark}

It is easy to check that  Lemma 2.1 and Lemma 2.2 in \cite{LT05}
still hold for the twisted Courant bracket, and that the same proofs
work. More specifically, if $f: M \to g^*$ is a submersion, then
$df_\C^\perp \subset T_\C M \oplus T^*_\C M$ is closed under the
$H$-twisted Courant bracket for any closed $H \in \Omega^3(M)$.
Therefore, if $\Gamma$ is a sub-bundle of $T_\C M \oplus T^*_\C M$
which is closed under the $H$-twisted Courant bracket, then the
image of $\Gamma \cap df^\perp_\C$  in $T_\C (f^{-1}(0)) \oplus
T^*_\C(f^{-1}(0))$ is closed under the $H|_{f^{-1}(0)}$-twisted
Courant bracket. Similarly, let $G$ act freely on $M$ and assume
that $H \in \Omega^3(M)$ is closed and basic. Then $H$ descends to a
closed form $\widetilde{H} \in \Omega^3(M/G)$ and the set of
$G$-invariant sections of $(\g_M)^\perp_\C$ is closed under the
$H$-twisted Courant bracket. Therefore, if $\Gamma$ is a sub-bundle
of $T_\C M \oplus T^*_\C M$ which is closed under the $H$-twisted
Courant bracket, then the image of $\Gamma \cap (\g_M)^\perp_\C$  in
$T_\C (M/G) \oplus T^*_\C(M/G)$ is  closed under the
$\widetilde{H}$-twisted Courant bracket.

In this context, we will work with a variant of the notion of
generalized moment map we defined in \cite{LT05}.

\begin{definition} \label{deftmm}
Let a compact Lie group $G$ with Lie algebra $\g$ act on a manifold
$M$, preserving an $H$-twisted generalized complex structure
$\mathcal{J}$, where $H \in \Omega^3(M)^G$ is closed. Let $L \subset
T_\C M \oplus T^*_\C M$ denote the $\sqrt{-1}$ eigenbundle of $\J$.
A {\bf twisted generalized moment map}  is a smooth function $f
\colon  M \to \g^*$ so that
\begin{itemize}
\item
There exists a one form $\alpha \in \Omega^1(M,\g^*)$, called the
{\bf moment one form}, so that $\xi_M - \newi (df^\xi+  \newi
\alpha^\xi)$ lies in $L$ for all $\xi \in \g$, where $\xi_M$ denote
the induced vector field.
\item $f$ is equivariant.
\end{itemize}
\end{definition}

\begin{example}  \
\begin{itemize}
\item [a)]
Let $G$ act on a generalized complex manifold with generalized
moment map $f + \newi h$. Then $f$ is a twisted generalized moment
map with moment one form $dh$.
\item [b)]
Let  $G$ act on an $H$-twisted generalized complex manifold $(M,\J)$
with twisted generalized moment map $f$ and moment one-form
$\alpha$. If $B \in \Omega^2(M)^G$, then $G$ acts on the
$B$-transform of $\J$ with twisted generalized moment map $f$ and
moment one form $\alpha'$, where $(\alpha')^\xi = \alpha^\xi +
\iota_{\xi_M} B$ for all $\xi \in \g$.
\end{itemize}
\end{example}

Let a compact Lie group $G$ act on a twisted generalized complex
manifold $(M,\J)$ with twisted generalized moment map $f$. Let
$\cO_a$ be the co-adjoint orbit through $a \in \g^*$. If $G$ acts
freely on $f^{-1}(\cO_a)$, then $\cO_a$ consists of regular values
and $M_a = f^{-1}(\cO_a)/G$ is a manifold, which we still call the
{\bf generalized complex quotient}. The proof of Lemma 3.8 in
\cite{LT05} applies almost word for word to the following
generalization.

\begin{lemma} \label{twisted complex reduction}
Let a compact Lie group $G$ act on an $H$-twisted generalized
complex manifold $(M,\J)$ with a  twisted generalized moment map $f
\colon M \to \g^*$. Let $\cO_a$ be a co-adjoint orbit through $a \in
\g^*$ so that $G$ acts freely on $f^{-1}(\cO_a)$. Assume that the
moment one-form is trivial and $H$ is basic. Then $H$  descends to
$\widetilde{H} \in \Omega^3_0(M_a)$ and  the generalized complex
quotient  $M_a$ naturally inherits a $\widetilde{H}$-twisted
generalized complex structure $\widetilde{\J}$.

Moreover, for all $m \in f^{-1}(\cO_a)$,
$$\type(\widetilde{\mathcal{J}})_{[m]} = \type(\mathcal{J})_m.$$
\end{lemma}

\begin{example} \
\itemize
\item [a)]
Let a compact Lie group  act on a (untwisted) generalized complex
manifold with a real generalized moment map. Then Lemma 3.8 in
\cite{LT05} and Lemma~\ref{twisted complex reduction} yield
identical generalized complex structures on $M_a$.
\item
[b)] Let a compact Lie group $G$ act on an $H$-twisted generalized
complex manifold $(M,\J)$ with twisted generalized moment map $f: M
\to \g^*$. Assume that the moment one-form is trivial and $H$ is
basic. If $B \in  \Omega^2(M)$ is basic, then it descends to
$\widetilde{B} \in \Omega^2(M_a)$, the $B$-transform  of $\J$
satisfies the same conditions, and its twisted generalized complex
quotient is the $\widetilde{B}$-transform of the twisted generalized
complex quotient of $\J$.
\end{example}

\begin{lemma}\label{twist three form}
Let a compact Lie group $G$ act freely on a manifold $M$. Let $H$ be
an invariant closed three form and let $\alpha$ be an equivariant
mapping from $\g$ to $\Omega^1(M)$.  Fix a connection $\theta \in
\Omega(M,\g)$. Then if $H + \alpha \in \Omega^3_G(M)$ is
equivariantly closed, there exists a natural form $\Gamma \in
\Omega^2(M)^G$ so that $\iota_\xi \Gamma = \alpha^\xi$. Thus $H +
\alpha + d_G \Gamma \in \Omega^3(M)^G \subset \Omega_G^3(M)$ is
closed and basic and so descends to a closed form $\widetilde{H} \in
\Omega^3(M/G)$ so that $[\widetilde{H}]$ is the image of $[H  +
\alpha]$ under the Kirwan map.\footnote{ In fact, the analogous
statement holds for any equivariantly closed form, but the proof in
general is more involved.}
\end{lemma}

\begin{proof}
Since $H + \alpha$ is equivariantly closed
\begin{itemize}
\item $dH = 0$.
\item $\iota_{\xi_M} H = d \alpha^\xi$ for all $\xi \in \g$.
\item $\iota_{\xi_M} \alpha^\eta = - \iota_{\eta_M} \alpha^\xi$ for
all $\xi$ and $\eta$ in $\g$.
\end{itemize}
Define $\beta \in \Omega^2(M)^G$ by $\beta(X,Y) = -
\alpha^{\theta(X)}(\theta(Y)_M)= \alpha^{\theta(Y)}(\theta(X)_M)$
for every vector field $X$ and $Y$. Define $\Gamma =  \Gamma_\theta
= -(\alpha,\theta) + \beta \in \Omega^2(M)^G$. Then $\iota_{\xi_M}
\Gamma = \alpha^\xi$ for all $\xi \in \g$. Since $\Gamma$ is
$G$-invariant, $\iota_{\xi_M} d \Gamma = -d \iota_{\xi_M} \Gamma = -
d \alpha^\xi = -\iota_{\xi_M}H$. Therefore, $H + d \Gamma$ is basic.
\end{proof}

We can now prove the twisted version of Proposition 3.10.
\begin{proposition} \label{Twisted Complex Reduction}
Let a compact Lie group $G$ act on an $H$-twisted generalized
complex manifold $(M,\mathcal{J})$ with twisted generalized moment
map $f: M \to \g^*$ and moment one-form $\alpha \in
\Omega^1(M,\g^*)$. Let $\cO_a$ be a co-adjoint orbit through $a \in
\g^*$ so that $G$ acts freely on $f^{-1}(\cO_a)$. Assume that $H +
\alpha$ is equivariantly closed. Given a connection on
$f^{-1}(\cO_a)$, the twisted generalized complex quotient $M_a$
inherits a $\widetilde{H}$-twisted generalized complex structure
$\widetilde{\J}$, where $\widetilde{H}$ is defined as in the Lemma
above. Up to $B$-transform, $\widetilde{J}$ is independent of the
choice of connection. Finally, for all $m \in f^{-1}(\cO_a)$,
$$\type(\widetilde{\mathcal{J}})_{[m]} = \type(\mathcal{J})_m.$$
\end{proposition}

\begin{proof}

By restricting to a neighborhood of $f^{-1}(\cO_a)$, we may assume
that $G$ acts freely on $M$. Given a connection $\theta \in
\Omega^1(M,\g)$, by the lemma above there exists $\Gamma \in
\Omega^2(M)$ so that $\iota_{\xi_M} \Gamma = \alpha^\xi$. So the
$\Gamma$-transform of $\J$ is a $H - d \Gamma$-twisted generalized
complex structure with twisted generalized  moment map $f$ and
trivial moment one-form, and hence induces a natural
$\widetilde{H}$-twisted generalized complex structure
$\widetilde{\J}$ on $M_a$ by Lemma ~\ref{twisted complex reduction}.

Given another connection $\theta'$, then $\Gamma_{\theta} -
\Gamma_{\theta'}$ is basic and hence descends to $\widetilde{\gamma}
\in \Omega^2(M)$. Consequently, the resulting twisted generalized
complex structure on $M_a$ is the $\widetilde{\gamma}$ transform of
$\widetilde{\J}$.
\end{proof}

\begin{remark}
In particular, even if $H = 0$, if $[\alpha] \neq 0$ then in general
$[\widetilde{H}]$ will also not vanish. Thus, in principle it may be
possible to get non-trivially twisted quotients from non-twisted
spaces.
\end{remark}

\begin{example} \
\begin{itemize}
\item [a)]
Let a compact Lie group act on a generalized complex manifold with
generalized moment map $f + \newi h$. The generalized complex
structure on $M_a$ induced by Proposition~\ref{Twisted Complex
Reduction} is the $ (h,d \theta)$-transform of the generalized
complex structure induced by Proposition~\ref{Twisted Complex
Reduction}. In particular, it may be twisted.
\item [b)]
Let $G$ act on an $H$-twisted generalized complex manifold $(M,\J)$
with twisted generalized moment map $f$ and moment one-form
$\alpha$. Assume  that $H + \alpha$ is equivariantly closed. Fix a
connection $\theta$ on $f^{-1}(\cO_a)$.

Given any $B \in \Omega^2(M)^G$, let $\Gamma \in \Omega^2(M)^G$ be
the natural form associated to the closed form $d_G(B) \in
\Omega^3_G(M)$ by  Lemma~\ref{twist three form} . Then $B + \Gamma$
is basic and hence descends to a form $\widetilde{B} \in
\Omega^2(M_a)$, and the twisted generalized complex quotient of the
$B$-transform of $\J$ is the $\widetilde{B}$ transform of the
twisted generalized complex quotient. Hence,  we cannot get
interesting new examples by  applying $B$-transforms to the space
upstairs.
\end{itemize}
\end{example}

A {\bf twisted generalized moment map} and {\bf moment one-form} for
a group action on a twisted generalized K\"ahler manifold
$(M,\J_1,\J_2)$ are simply a twisted generalized moment map and
moment one-form for the twisted generalized complex structure
$\J_1$. We are now ready to prove our final proposition.

\begin{proposition} \label{tgkahler quotient}
Let a compact connected Lie group $G$ act on an $H$-twisted
generalized K\"{a}hler manifold $(M,\J_1,\J_2)$ with twisted
generalized moment map $f \colon M \to \g^*$  and moment one-form
$\alpha \in \Omega^1(M,\g^*)$. Let $\cO_a$ be a co-adjoint orbit
through $a \in \g^*$ so that  $G$ acts freely on $f^{-1}(\cO_a)$.
Assume $H+ \alpha \in \Omega_G^3(M)$ is equivariantly closed. Then
the generalized K\"ahler quotient $M_a$ naturally inherits a
$\widetilde{H}$-twisted generalized K\"{a}hler structure
$(\widetilde{\J}_1,\widetilde{\J}_2)$, where $\widetilde{H}$ is
defined as in Lemma~\ref{twist three form} using the canonical
connection on $M$.

Finally, let $\mathfrak{k}$ be the Lie algebra of the stabilizer $K$
of $a$, and let $L_2$ be the $\sqrt{-1}$ eigenbundle of $\J_2$. Then
for all $m \in M$,
$$\type(\widetilde{\J}_1)_{[m]} = \type(\J_1)_m, \qquad \mbox{and}$$
$$\type(\widetilde{\J}_2)_{[m]} = \type(\J_2)_m - \frac{1}{2}\dim G
- \frac{1}{2}\dim(K) + 2 \dim (\mathfrak{h}_M \cap \pi(L_2))_m.$$
\end{proposition}

\begin{proof}
By restricting to a neighborhood of $f^{-1}(\cO_a)$, we may assume
that $G$ acts freely on $M$. The $\Gamma_{\theta}$ transform of
$(\J_1,\J_2)$ is a $H +d\Gamma_{\theta}$-twisted generalized
K\"ahler structure with twisted generalized moment map $f$ and
trivial moment one-form. From this point, the proof of Proposition
4.6 in \cite{LT05} goes through without change.
\end{proof}

\end{document}